\documentclass[12pt]{amsart}

\newtheorem{theorem}{Theorem} [section]
\newtheorem{prop}[theorem]{Proposition}
\newtheorem{lemma}[theorem]{Lemma}
\newtheorem{cor}[theorem]{Corollary}

\newtheorem{question}[theorem]{Question}

\numberwithin{equation}{section}

\newcommand{\E}{\mathbb{E}}
\newcommand\A{{\mathbb A}}
\newcommand\C{{\mathbb C}}

\renewcommand\P{{\mathbb P}}
\newcommand\R{{\mathbb R}}
\newcommand\Z{{\mathbb Z}}
\newcommand\Q{{\mathbb Q}}
\newcommand\F{{\mathbb F}}

\newcommand\tQ{\widetilde{\mathbb Q}}
\newcommand\tF{\widetilde{\mathbb F}}
\newcommand\hO{\widehat{\bf O}}
\newcommand\hF{\widehat{F}}
\newcommand\cV{{\mathcal V}}
\newcommand\cK{{\mathcal K}}
\newcommand\cG{{\mathcal G}}

\newcommand\eps{\varepsilon}

\renewcommand\phi{\varphi}
\newcommand\iso{\simeq}

\newcommand\Gal{\operatorname{Gal}}
\newcommand\vol {\operatorname{vol}}
\newcommand\Res {\operatorname{Res}} 
\newcommand\Spec{\operatorname{Spec}}
\newcommand\Det{\operatorname{Det}}

\begin{document}

\thispagestyle{empty}

\title[Transfinite diameter and the resultant]
{Transfinite diameter and the resultant}

\author{Laura DeMarco}
\address{Laura DeMarco\\
Department of Mathematics\\
University of Chicago\\
5734 S. University Avenue\\
Chicago, IL 60637\\
USA}
\email{demarco@math.uchicago.edu}

\author{Robert Rumely}
\address{Robert Rumely\\ 
Department of Mathematics\\
University of Georgia\\
Athens, GA 30602\\
USA}
\email{rr@math.uga.edu}

\date{January 4, 2006}
\thanks{The first author was supported in part by an NSF Postdoctoral Fellowship
(DMS 0303421)
and the second author by NSF grant DMS 0300784.}

\begin{abstract} 
We prove a formula for the Fekete-Leja transfinite diameter
of the pullback of a set $E \subset \C^N$ by a regular polynomial map $F$, 
expressing it in terms of the resultant of  the leading part of $F$ and the
transfinite diameter of $E$.  We also establish the nonarchimedean analogue
of this formula.  A key step in the proof is a formula for the 
transfinite diameter of the filled Julia set of $F$.  
\end{abstract} 

\maketitle

\section{Introduction}
A polynomial map $F:\C^N\to\C^N$ is {\em regular} if it 
extends to a holomorphic endomorphism of $\P^N$.  In that case, 
its coordinate functions necessarily all have the same degree, and 
the {\em degree} of $F$ is their common degree.  
In this note, we establish the following pullback formula for the 
Fekete-Leja transfinite diameter:  

\begin{theorem}  \label{pullback}
For any regular polynomial map 
$F:\C^N\to\C^N$ of degree $d$, and any bounded set $E \subset \C^N$, 
 $$d_\infty(F^{-1}E) = |\Res(F_h)|^{-1/Nd^N} d_\infty(E)^{1/d} \ .$$ 
\end{theorem}

\noindent  
Here, $d_{\infty}(E)$ is the Fekete-Leja transfinite diameter of $E$ 
(see \S\ref{capacity}),
and $\Res(F_h)$ is the multiresultant (\cite{vdW}, \cite{GKZ}) 
of the $N$ homogeneous polynomials comprising the leading part of $F$ ;
it is a homogeneous polynomial of degree $Nd^{N-1}$ in the coefficients of $F_h$ 
which is nonzero precisely when $F$ is regular (see \S\ref{resultant}).   

Theorem \ref{pullback} generalizes a classical 
result in dimension 1, due to Fekete, where the resultant is the leading
coefficient of the polynomial $F$ (see \cite{Goluzin}).  When 
$F$ is a linear automorphism of $\C^N$, its resultant is the 
determinant, and the formula is due to She\u{i}nov \cite{Sheinov}.  
In dimension $N=2$, in the case 
when $F$ is defined over a number field and $E$ is a polydisc, 
it was shown by Baker and Rumely \cite{Baker:Rumely:equidistribution}.
Other special cases were proved by Bloom and Calvi \cite{Bloom:Calvi}, 
who conjectured the existence of a general pullback formula.
\vskip .1 in

For each prime $p$, let $\C_p$ be the field of $p$-adic complex numbers,
the completion of the algebraic closure of 
the field of $p$-adic numbers $\Q_p$, 
equipped with its nonarchimedean absolute value $|x|_p$.  
There is a natural analogue of the Fekete-Leja transfinite diameter 
for sets in $\C_p^N$ (see \S\ref{vbackground}).  Using Theorem \ref{pullback}, 
together with an approximation theorem of Moret-Bailly \cite{MB},
we deduce the nonarchimedean counterpart to Theorem \ref{pullback}: 

\begin{theorem} \label{pullbackv}
For any regular polynomial map 
$F:\C_p^N\to\C_p^N$ of degree $d$, and any bounded set $E \subset \C_p^N$, 
 $$d_\infty(F^{-1}E)_p = |\Res(F_h)|_p^{-1/Nd^N} d_\infty(E)_p^{1/d} \ .$$ 
\end{theorem}

\bigskip
Theorem \ref{pullback} also yields explicit formulas for the transfinite diameter of 
certain sets in $\C^N$.  Here we give two examples, extending results 
in \cite{Bloom:Calvi}:   

\begin{cor} \label{AnPoly} {\em (Special analytic polyhedra)}
Suppose $F = (F_1, \ldots , F_N)$ is a regular polynomial map of degree 
$d$, and let $K = \{z\in\C^N: |F_i(z)| \leq 1 \mbox{ for all } i\}$.  Then 
  $$ d_\infty(K) = |\Res(F_h)|^{-1/Nd^N} \ .$$
\end{cor}

\proof
The unit polydisc $D(0,1)$ has transfinite diameter 1, and $K = F^{-1}D(0,1)$.
\qed

\begin{cor}  \label{filledJ}  {\em (Filled Julia sets)}
Suppose $F$ is a regular polynomial map of degree $d \geq 2$, and let 
$K_F = \{z\in\C^N: \sup_n \|F^n(z)\| < \infty \}$ be its filled Julia set.  Then
   $$d_\infty(K_F) = |\Res(F_h)|^{-1/Nd^{N-1}(d-1)}.$$
\end{cor}

\proof
The filled Julia set is compact, non-pluripolar,   
and satisfies $F^{-1} K_F = K_F$.
\qed

\medskip
When $N=1$, Corollary \ref{filledJ} is well-known (see e.g. \cite{Ransford}).  
In dimension $N=2$, when $F$ is homogeneous, it was proved with 
complex analytic methods by DeMarco for the homogeneous capacity of $K_F$, 
and used to study the bifurcation current for holomorphic families 
of rational maps \cite{D:lyap}. Baker and Rumely later showed that in $\C^2$,  
DeMarco's homogeneous capacity coincides with the transfinite diameter 
\cite{Baker:Rumely:equidistribution}.

\bigskip\noindent{\bf Outline of the proof of Theorem \ref{pullback}.} 
The proof uses a combination of analytic and arithmetic techniques.  
We first prove Corollary \ref{filledJ} for homogeneous 
regular polynomial maps by combining two recent results, 
one in dynamics and the other in number theory.  
In a study of the Lyapunov exponents of a homogeneous polynomial 
map $F$, Bassanelli and Berteloot \cite{Bassanelli:Berteloot} 
equated the resultant of $F$ 
to a potential-theoretic expression.     
Using arithmetic intersection theory, 
Rumely \cite{Rumely:Vformula} gave a formula for the transfinite
diameter of a compact set $E \subset \C^N$ in terms of its pluricomplex Green's 
function. When $E = K_F$ these formulas can be related by  
integration by parts, yielding the formula for $d_\infty(K_F)$ 
(Theorem \ref{Julia}).  

Next, using the pullback formula for the global sectional capacity 
(\cite{Rumely:Lau:Varley}, Theorem 10.1),
we show that when $F$ is homogenous and 
defined over the number field $\Q(i)$, there is a 
constant $C_F$ such that for any compact set $E \subset \C^N$ 
$$ d_{\infty}(F^{-1}E) = C_F \cdot d_{\infty}(E)^{1/d}.$$  
Since the filled Julia set satisfies $F^{-1}K_F = K_F$, 
Theorem \ref{Julia} lets us evaluate $C_F = |\Res(F)|^{-1/Nd^{N}}$.  
Finally an approximation argument, using the continuity of the resultant
and properties of the transfinite diameter, gives Theorem \ref{pullback}.


\bigskip
\section{The transfinite diameter in $\C^N$}
\label{capacity}

\noindent{\bf The transfinite diameter.}
For each fixed integer $n\geq0$, we let $e_1, e_2, \ldots, e_{M(n)}$ denote
the monomials in $N$ variables of degree 
$\leq n$, in any order.   Note that 
 $$M(n) = \left(\!\!\begin{array}{c} N+n \\ N \end{array} \!\!\right)$$  
coincides with the dimension of the space of homogeneous polynomials of degree $n$ 
in $N+1$ variables.  Given a bounded set $E\subset\C^N$, 
its $n$-th  diameter is defined to be: 
\begin{equation}  \label{Vand}
  d_n(E) = \left(\sup_{\zeta_1, \ldots, \zeta_{M(n)}\in E} 
   | \Det(e_i(\zeta_j))_{i,j}|\right)^{1/D(n)},
\end{equation}
where
  $$D(n) = 
    \sum_{m=1}^{n} m\left(\!\!\begin{array}{c} N+m-1 \\ m \end{array}\!\!\right) 
     = N \! \left(\!\!\begin{array}{c} N+n \\ N+1 \end{array}\!\!\right)$$
is the degree of the Vandermonde determinant in equation (\ref{Vand}).  
The {\em transfinite diameter} of $E$ is the limit 
\begin{equation} \label{diam}
  d_\infty(E) = \lim_{n\to\infty} d_n(E).
\end{equation}
The  existence of the limit in (\ref{diam}) 
was posed as a question by Leja in 1959, and was established  
for every bounded set by Zaharjuta \cite{Zaharjuta} in 1975.  

\bigskip\noindent{\bf Dimension $N=1$.}
In dimension one, the definitions reduce to the familiar transfinite diameter
in $\C$; namely, 
  $$d_n(E) = \left(\sup \prod_{i<j} |\zeta_i-\zeta_j|\right)^{2/n(n+1)},$$
where the supremum is taken over all sets of $n+1$ points in $E$, and 
  $$d_\infty(E) = \lim_{n\to\infty} d_n(E).$$
Recall that in dimension one, if $E$ is compact, 
the transfinite diameter coincides with the 
logarithmic capacity, which can be computed in terms of the Robin 
constant (see \cite{Ahlfors:conformal}).  Namely, if $E\subset\C$ is compact, 
there is a constant $V(E)$ such that the Green's function satisfies  
  $$G_E(z) = \log|z| + V(E) + o(1)$$
for $z$ near $\infty$, and then 
\begin{equation} \label{Robin1}
  d_\infty(E) = e^{-V(E)} \ .
\end{equation}

\bigskip\noindent {\bf The higher dimensional Robin formula.}
The Robin formula (\ref{Robin1}) for the transfinite diameter 
was generalized to arbitrary dimensions in \cite{Rumely:Vformula}.
Let $E$ be a compact set in $\C^N$, and let $G_E$ 
be its pluricomplex Green's function.  
That is, with $\|z\| = \sqrt{|z_1|^2+ \cdots + |z_N|^2}$, 
      $$G_E(z) = 
  	\left( \sup\{u(z): u\in PSH(\C^N), u \leq \log^+\|\cdot\| + O(1) \mbox{ and } u|_E\leq 0\}
	\right)^*.$$
By definition,  $E$ is 
non-pluripolar if $G_E\not\equiv \infty$;  
this is equivalent to $d_{\infty}(E) > 0$. (See \cite{Klimek} for general
properties of the pluricomplex Green's function.)

Identify
$\C^N$ with the affine chart $\A_0\subset \P^N$, so that  
$(z_1, \ldots, z_N) = (Z_0: \ldots: Z_N)$ are the coordinates on 
$\C^N$ when $Z_0\not=0$.  The hyperplane $H_0 = \{Z_0=0\}$ in 
$\P^N$ will be identified with $\P^{N-1}$.  For each $j = 1, \ldots, N$,
define functions on $\P^{N-1}$ by  
 $$g_j(z_1:\cdots:z_N) = \limsup_{|t|\to\infty} \,(G_E(tz) - \log|tz_j|)$$
and set
  $$g_E(z_1:\cdots:z_N) = \limsup_{|t|\to\infty}\, (G_E(tz) -\log\|tz\|).$$
Then $dd^c g_j = T_E - T_j$ where $T_E$ is a well-defined positive 
$(1,1)$-current on $\P^{N-1}$ and $T_j$ 
is the current of integration over the hyperplane $H_j = \{z_j=0\}$;  
and $dd^c g_E = T_E - \omega$, where $\omega$ is the Fubini-Study form.  

Rumely (\cite{Rumely:Vformula}, Theorem 0.1) showed that if $E\subset\C^N$
is compact and non-pluripolar, then
\begin{equation} \label{R}
  - \log d_\infty(E) = \frac{1}{N} \sum_{j=0}^{N-1} \int_{\P^{N-1}} 
   g_{N-j} \,\, T_E^j\wedge T_1\wedge\cdots\wedge T_{N-j-1}.
\end{equation}
Observe that $T_1\wedge \cdots \wedge T_{N-j-1}$ is the current of 
integration along the $j$-dimensional linear subspace $\{z_1 = \cdots
= z_{N-j-1} = 0\}$ in $\P^{N-1}$.

\begin{question}
Is there a direct complex-analytic proof of equation (\ref{R})?
\end{question}


\bigskip
\section{Formulas for the transfinite diameter}
\label{Robin}

In this section we simplify the generalized Robin formula (\ref{R}).

\begin{prop} \label{Vformula}
For any compact, non-pluripolar set $E\subset\C^N$, the transfinite diameter
$d_\infty(E)$ satisfies: 
\begin{equation*} 
-\log d_\infty(E) 
 =  \frac{1}{N} \sum_{j=0}^{N-1}\int_{\P^{N-1}} g_E \,\, \omega^j \wedge T_E^{N-j-1}
	\,\, + \, \, \frac12 \sum_{j=2}^N \frac{1}{j} \ . 
\end{equation*}    
\end{prop}

\medskip
\proof
The desired expression follows from (\ref{R})
and integration by parts.  We will use induction to show that 
\begin{equation} \label{stepn}
-\log d_\infty(E) = \frac{1}{N}\sum_{j=0}^{n-1} \int_{\P^{N-1}} g_E \,\, 
                     \omega^j\wedge T_E^{N-1-j}  \qquad\qquad 
\end{equation}
 $$	\qquad\qquad\qquad + \,\,\, \frac{1}{N}\sum_{j=0}^{n-1} \int_{\P^{N-1}} (g_{N-j}-g_E)
		     \,\, \omega^j \wedge T_1\wedge\cdots\wedge T_{N-1-j} $$
 $$  \qquad\qquad\qquad\quad + \,\,\, \frac{1}{N} \sum_{j=0}^{N-n-1} \int_{\P^{N-1}} g_{j+1} \, \,
		     \omega^n\wedge T_E^{N-1-j-n} \wedge T_1 \wedge\cdots
		     \wedge T_j, $$
for each $n = 0, \ldots, N$.
The base case $n=0$ is equation (\ref{R}).  

To pass from $n$ to $n+1$, we use the identity 
 $$\omega^n\wedge T_E^{N-n-1} - \omega^n\wedge T_1\wedge\cdots\wedge T_{N-n-1} 
   \quad\quad\quad $$  $$ \quad\quad\quad = \sum_{j=0}^{N-n-2} dd^c g_{j+1} 
\wedge\omega^n \wedge T_E^{N-n-j-2}\wedge T_1\wedge\cdots\wedge T_j.$$
Integrating the identity against $g_E$, we can apply integration-by-parts 
to move the $dd^c$ from the $g_{j+1}$ to the $g_E$.  This yields, 
\begin{equation} \label{IBP}
 \sum_{j=0}^{N-n-2} \int_{\P^{N-1}} g_{j+1} \,\, (T_E - \omega)
\wedge\omega^n \wedge T_E^{N-n-j-2}\wedge T_1\wedge\cdots\wedge T_j \qquad\qquad 
\end{equation}
$$ \qquad\qquad  = \int_{\P^{N-1}} g_E \,\, \omega^n\wedge T_E^{N-n-1}
   - \int_{\P^{N-1}} g_E \,\, \omega^n\wedge T_1\wedge\cdots\wedge T_{N-n-1}.$$
To justify the integration by parts, we refer the reader to \cite{Bedford:Taylor}, \S2,
and \cite{Demailly:lelong}, \S3.
Note that the function $-g_{j+1}$ on $\P^{N-1}$ 
(which satisfies $dd^c (-g_{j+1}) = T_{j+1}-T_E$) 
can be approximated by a decreasing sequence $\{g_{j+1,k}\}_{k\geq 0}$ of 
bounded functions such that 
  $$dd^c g_{j+1,k} + T_E \geq 0$$
for all $k$.  
The operator $dd^c$ can be moved from $g_{j+1,k}$ to $g_E$ 
and then we can pass to the limit.  

Substituting equation (\ref{IBP}) into the final sum of equation (\ref{stepn}) for $n$, 
we obtain equation (\ref{stepn}) for $n+1$.  

For $n=N$, the inductive expression (\ref{stepn}) gives us
\begin{eqnarray*}
-\log d_\infty(E) &=& \frac{1}{N} \sum_{j=0}^{N-1} \int_{\P^{N-1}} 
       g_E \,\, \omega^j \wedge T_E^{N-1-j} \\
&+& \frac{1}{N} \sum_{j=0}^{N-1} \int_{\P^{N-1}} (g_{N-j} - g_E) \,\,
    \omega^j\wedge T_1\wedge\cdots\wedge T_{N-j-1}.
\end{eqnarray*}
The final summation is independent of the set $E$ and can be
rewritten as 
 $$\frac{1}{N}\sum_{j=0}^{N-1} \int_{\P^{N-1}} (\log\|z\| - \log|z_{N-j}|) \,\,
    \omega^j\wedge T_1\wedge\cdots\wedge T_{N-j-1} \qquad\qquad\qquad$$
\begin{eqnarray*}
	&=& \frac{1}{N} \sum_{j=1}^{N-1} \int_{\P^j} (\log\|z\|-\log|z_1|) \,\,\omega^j  
   \ = \  - \frac{1}{N} \sum_{j=1}^{N-1} \int_{\C^{j+1}} \log|z_1| \,\, dm \\
   &=& \frac{1}{2N} \sum_{j=1}^{N-1} \left( 1 + \frac12 + \cdots + \frac{1}{j} \right) 
   \ = \ \frac{1}{2N} \sum_{j=1}^{N-1} \frac{N-j}{j} 
   \ = \ \frac{1}{2} \sum_{j=2}^N \frac{1}{j}  \ , 
\end{eqnarray*}
where $dm$ on the second line is the standard area form on the unit sphere of unit volume.
\qed

\bigskip 
\noindent {\bf A symmetric form of the Robin formula.}  
Let $B^N = \{z \in \C^N : \|z\| \le 1\}$ be the $L^2$-unit ball in $\C^N$.  Then 
$G_{B^N}(z) = \log^{+}(\|z\|)$, $g_{B^N}(z) = 0$, and $T_{B^N} = \omega$.  
Taking $E = B^N$
in Proposition \ref{Vformula}, we recover the following fact 
(\cite{Jedrz}, \cite[p.555]{Rumely:Lau}): 

\begin{cor} \label{BNdiam} \quad 
$d_{\infty}(B^N) \ = \ \exp \left(- \frac12 \sum_{j=2}^N \frac{1}{j} \right).$
\end{cor}  

\medskip\noindent
Thus, Proposition \ref{Vformula} can be reformulated as saying that 
for an arbitrary compact, non-pluripolar set $E\subset\C^N$, 
if  $\omega_E = \frac{1}{N} \sum_{j=0}^{N-1}  T_E^{N-j-1} \wedge \omega^j$, then 
\begin{equation} \label{VformulaII} 
-\log(d_{\infty}(E)) 
\ = \ -\log(d_{\infty}(B^N)) \ + \ \int_{\P^{N-1}} g_E \, \omega_E \ .     
\end{equation}

This can be generalized as follows:

\begin{theorem} \label{GenVFormula} 
Fix a compact, non-pluripolar set $E_0 \subset \C^N$.  
For an arbitrary compact non-pluripolar set $E \subset \C^N$, 
put $\omega_{E,E_0} = \frac{1}{N} \sum_{j=0}^{N-1}  T_E^{N-j-1} \wedge T_{E_0}^j$. 
Then 
\begin{equation} \label{SymVF} 
-\log(d_{\infty}(E)) \ = \ -\log(d_{\infty}(E_0)) 
                     \ + \ \int_{\P^{N-1}} (g_E\!-\!g_{E_0}) \, \omega_{E,E_0} \ .  
\end{equation}                      
\end{theorem} 

\smallskip
\proof
By repeatedly using the identity 
$\omega = T_{E_0} - (T_{E_0} - \omega) = T_{E_0} - dd^c g_{E_0}$, one gets 
\begin{equation} \label{EFormula} 
g_E \omega_E 
 \ = \ \frac{1}{N} \sum_{j=0}^{N-1} g_E \, T_E^{N-j-1} \wedge T_{E_0}^{j} 
           - \frac{1}{N} \sum_{i+j+\ell = N-2} 
         g_E \, dd^c g_{E_0} \wedge T_E^{i} \wedge T_{E_0}^j \wedge \omega^{\ell} \ . 
\end{equation}
Similarly, using $\omega = T_E - dd^c g_E$, one finds that    
\begin{equation} \label{E0Formula} 
g_{E_0} \omega_{E_0} 
 \ = \ \frac{1}{N} \sum_{j=0}^{N-1} g_{E_0} \, T_E^{N-j-1} \wedge T_{E_0}^j
              - \frac{1}{N} \sum_{i+j+\ell = N-2} 
         g_{E_0} \, dd^c g_E \wedge T_E^{i} \wedge T_{E_0}^j \wedge \omega^{\ell} \ . 
\end{equation}
Inserting (\ref{EFormula}) and (\ref{E0Formula}) into the expressions (\ref{VformulaII}) 
corresponding to $E$ and $E_0$, then subtracting and using integration by parts, one
obtains (\ref{SymVF}).  The integration by parts is justified by the same 
remarks as in the proof of Proposition \ref{Vformula}.  
\qed


\bigskip
\section{Regular polynomial endomorphisms and the resultant}
\label{resultant}

A polynomial map $F:\C^N\to\C^N$ is {\em regular} if it 
extends to a holomorphic endomorphism of $\P^N$.  
This map is necessarily finite.  The {\em degree} of 
$F$ is the degree of each of its coordinate functions.  
The polynomial map $F:\C^N\to\C^N$ is regular of degree $d \geq 1$
if and only if each component has degree $d$ and the leading homogeneous part
$F_h$ satisfies $F_h^{-1}\{0\} = \{0\}$.  
Alternatively, 
  $$\liminf_{\|z\|\to\infty} \frac{\|F(z)\|}{\|z\|^d} > 0.$$
Write $\log^+(x) = \max(0,\log(x))$, where $\log(x)$ is the natural logarithm.  
When $d \geq 2$ the {\em escape-rate function} 
$G^F: \C^N\to \R_{\geq 0}$ is defined by 
  $$G^F(z) = \lim_{n\to\infty} \frac{1}{d^n} \log^+ \|F^n(z)\|,$$
and $G^F$ coincides with the pluricomplex Green's function $G_{K_F}$ for the filled 
Julia set $K_F$ (see \cite{Fornaess:Sibony} and \cite{Bedford:Jonsson}).  
The escape-rate
function is continuous on $\C^N$ and maximally plurisubharmonic in the complement 
of $K_F$.  

\bigskip\noindent{\bf The resultant of a homogeneous polynomial map.}
When $F = (F_1, \ldots, F_N)$ is homogeneous, one defines its 
{\em resultant} $\Res(F)$ to be the multiresultant of its $N$ coordinate functions
(\cite{GKZ}, \cite{vdW}), namely the unique polynomial 
in the coefficients of $F$ such that    
\begin{itemize}
\item[(i)]  $\Res(F) \not= 0$ if and only if $F$ is regular, and 
\item[(ii)] $\Res (z_1^d, \ldots, z_N^d) = 1$.  
\end{itemize}
The polynomial $\Res$ is absolutely irreducible, homogeneous of degree $d^{N-1}$ 
in the coefficients of each $F_i$, 
and has total degree $Nd^{N-1}$ (\cite{GKZ}, Chapter 13).  

When $F$ is linear, $\Res(F) = \Det(F)$.  
When $F : \C^2 \rightarrow \C^2$ 
is a quadratic map with coordinate functions 
$F_1(z_1,z_2) = a_1 z_1^2 + b_1 z_1 z_2 + c_1 z_2^2$, 
$F_2(z_1,z_2) = a_2 z_1^2 + b_2 z_1 z_2 + c_2 z_2^2$, 
$$
\Res(F) = a_1^2 c_2^2 - 2 a_1 a_2 c_1 c_2 + a_2^2 c_1^2 
       - a_1 b_1 b_2 c_2 - a_2 b_1 b_2 c_1 + a_1 b_2^2 c_1 + a_2 b_1^2 c_2.
$$
For a quadratic map $F : \C^3 \rightarrow \C^3$, $\Res(F)$ is a homogeneous
polynomial of degree $12$ with $21894$ terms.  
 
\bigskip
\noindent{\bf The Bassanelli-Berteloot  formula.}
Let $F:\C^N\to\C^N$ be a regular homogeneous polynomial
map and $f:\P^{N-1}\to\P^{N-1}$ the induced holomorphic map 
on the hyperplane at infinity.  That is, $\pi\circ F = f\circ\pi$ where
$\pi:\C^N\setminus 0\to\P^{N-1}$ is the natural projection.  
On $\P^{N-1}$, define 
$$g^F(z_1:\cdots:z_N) = \limsup_{|t|\to\infty}\, (G^F(tz) -\log\|tz\|),$$
where $G^F$ is the escape rate function defined above.  
Then 
 $$dd^c g^F = T_f - \omega,$$
where $\omega$ is the Fubini-Study form on $\P^{N-1}$ and $T_f$
is a positive $(1,1)$-current, known as the {\em Green current} of $f$
(see \cite{Fornaess:Sibony}).

Bassanelli and Berteloot have shown (\cite{Bassanelli:Berteloot}, Proposition 4.9), 
\begin{equation} \label{BB}
  \sum_{j=0}^{N-1}\int_{\P^{N-1}} g^F \,\, \omega^j \wedge T_f^{N-j-1}
  = \frac{1}{d^{N-1}(d-1)} \log |\Res(F)| - \frac12 \sum_{j=1}^{N-1} \frac{N-j}{j} \ .
\end{equation}

\bigskip\noindent {\bf The filled Julia set.}
Using the Bassanelli-Berteloot formula (\ref{BB}), 
we can compute the transfinite diameter
of the filled Julia set of a homogeneous regular polynomial map
$F:\C^N\to\C^N$.  This special case of Corollary \ref{filledJ}
is the key to the proof of Theorem \ref{pullback}:  

\begin{theorem} \label{Julia}  
For each homogeneous regular polynomial 
map $F:\C^N\to\C^N$ of degree $d\geq 2$, 
the transfinite diameter of its filled Julia set is 
  $$d_\infty(K_F) = |\Res(F)|^{-1/Nd^{N-1}(d-1)}.$$
\end{theorem}

\proof 
For the filled Julia set $K_F$ of a regular homogeneous polynomial map $F$,
we have $g_{K_F} = g^F$ and $T_{K_F}=T_f$.  Thus 
Proposition \ref{Vformula} and the Bassanelli-Berteloot formula (\ref{BB}) yield
\begin{eqnarray*}
-\log d_\infty(K_F) &=& 
    \frac{1}{N} \sum_{j=0}^{N-1}\int_{\P^{N-1}} g^F \,\, \omega^j \wedge T_f^{N-j-1}
	 \, + \,  \frac{1}{2} \sum_{j=2}^N \frac{1}{j} \\ 
	 &=& \frac{1}{Nd^{N-1}(d-1)} \log |\Res(F)| \, - \, \frac{1}{2N} \sum_{j=1}^{N-1} \frac{N-j}{j} 
	 \, + \, \frac{1}{2} \sum_{j=2}^N \frac{1}{j} \\
    & =& \frac{1}{Nd^{N-1}(d-1)} \log |\Res(F)|.
\end{eqnarray*}
\qed


\bigskip
\section{Sectional capacity and the nonarchimedean transfinite diameter}
\label{vbackground}

In this section, we provide the arithmetic background needed for the proofs of 
Theorems \ref{pullback} and \ref{pullbackv}.  

For each rational prime $p$, let $\Q_p$ denote the field of $p$-adic numbers, 
equipped with its absolute value $|x|_p$ normalized so that $|p|_p = 1/p$.
Let $\tQ_p$ be the algebraic closure of $\Q_p$, and $\C_p$ its completion.
Then $\C_p$ is complete and algebraically closed.  The absolute value on $\Q_p$
extends in a unique way to $\C_p$, where it is still denoted $|x|_p$.  Write
$\hO_p$ for the ring of integers $\{x \in \C_p : |x|_p \leq 1\}$ 
and $m_p$ for its unique maximal ideal;  then $\hO_p/m_p \cong \tF_p$, the algebraic 
closure of the finite field $\F_p$.   
   
If $k$ is a number field, then for each place $v$ of $k$ there is a canonical
absolute value $|x|_v$ on the completion $k_v$, given by the modulus of additive
Haar measure.  
If $v$ is archimedean then $|x|_v = |x|^{[k_v:\R]}$;  if $v$ is
nonarchimedean, then $|x|_v = |x|_p^{[k_v:\Q_p]}$ if $v$ lies over the prime $p$ 
and $k_v$ is viewed as embedded in $\C_p$.  Each $|x|_v$ extends uniquely to an
absolute value on $\C_p$, and the relation $|x|_v = |x|_p^{[k_v:\Q_p]}$ continues 
to hold.  For each nonzero $\kappa \in k$, 
the {\em product formula} says that $$\prod_v |\kappa|_v = 1. $$ 
We will write $\C_v = \C_p$ if $v$ lies over $p$, 
and let $\Gal_c(\C_v/k_v) = \Gal_c(\C_p/k_v)$
be the group of continuous automorphisms fixing $k_v$.   
   
By definition, a local field (of characteristic $0$) 
is a finite extension of $\Q_p$ or $\R$.  
For any local field, there are infinitely many number fields $k$ with a place 
$v$ such that $k_v$ is isomorphic to the given local field.      
 
\bigskip\noindent {\bf Nonarchimedean transfinite diameter.}
Given a bounded set $E \subset \C_p^N$, one can define a transfinite diameter
$d_{\infty}(E)_p$ just as in the archimedean case, by replacing $|x|$ 
with $|x|_p$ in (\ref{Vand}):  for each $n$, put  
\begin{equation}  \label{Vandp}
 d_n(E)_p = \left(\sup_{\zeta_1, \ldots, \zeta_{M(n)} \in  E} 
| \Det(e_i(\zeta_j))_{i,j}|_p\right)^{1/D(n)},  
\end{equation}
where $M(n) = \left(\!\!\begin{array}{c} N+n \\ N \end{array}\!\!\right)$
and $D(n) = N \! \left(\!\!\begin{array}{c} N+n \\ N+1 \end{array}\!\!\right)$, 
and then let 
\begin{equation} \label{diamp} 
d_\infty(E)_p = \lim_{n\to\infty} d_n(E)_p\ .  
\end{equation}
The existence of the limit follows from (\cite{Rumely:Lau}, Theorem 2.6): in the
notation of \cite{Rumely:Lau}, if a local field $k_v$ is taken as the base field, 
and $H_0$ is the hyperplane at infinity for $\C_p^N\subset \P^N(\C_p)$,   
the limit  
\begin{equation} \label{DRel0}
d_{\infty}(E,H_0) := \lim_{n \rightarrow \infty}
      \left(\sup_{\zeta_1, \ldots, \zeta_{M(n)} \in  E} 
      | \Det(e_i(\zeta_j))_{i,j}|_v\right)^{(N+1)!/n^{N+1}}   
\end{equation}
exists; and since $\lim_{n \rightarrow \infty} D(n)\cdot (N+1)!/n^{N+1} = N$, 
while $|x|_v = |x|_p^{[k_v:\Q_p]}$, it follows that 
\begin{equation} \label{DRel} 
d_{\infty}(E)_p = d_{\infty}(E,H_0)^{1/(N [k_v:\Q_p])} \ .
\end{equation} 
Actually  \cite{Rumely:Lau} is written under a blanket hypothesis
that $E$ is stable under $\Gal_c(\C_p/k_v)$,    
but that assumption is not used in the proof of 
the existence of $d_{\infty}(E,H_0)$, 
which is a direct translation of Zaharjuta's proof over $\C$.  


\bigskip\noindent{\bf Approximation by polynomial polyhedra.}
Given a bounded set $E\subset \C_p^N$ and an $\varepsilon >0$, 
by \cite{Rumely:Lau}, Theorem 2.9, there is a finite
collection of polynomials $f_1, \ldots, f_M \in \C_p[z_1,\ldots,z_N]$ such that 
the set 
\begin{equation} \label{Nbhd} 
U = \{ z \in \C_p^N : |f_1(z)|_p \le 1, \ldots, |f_M(z)|_p \le 1 \} 
\end{equation} 
contains $E$ and satisfies 
 $$d_{\infty}(U)_p \le d_{\infty}(E)_p + \varepsilon \ ; $$
if $E$ is stable under $\Gal_c(\C_p/k_v)$ then one can take 
$f_1, \ldots, f_M \in k_v[z_1,\ldots,z_N]$.

\bigskip\noindent
{\bf Local sectional capacity.}
Rumely and Lau have shown that the transfinite diameter of a bounded set
$E\subset \C^N$ can also be computed as a growth rate of the set of polynomials
whose sup-norm is bounded by 1 on the set $E$, as the degree of the polynomials
tends to infinity.  More precisely, 
\begin{equation} \label{S}
  \log d_\infty(E) = 
    -\frac{1}{2N} \lim_{n\to\infty} \frac{(N+1)!}{n^{N+1}} \log \vol \{f\in\Gamma(n):  \|f\|_E\leq 1\},
\end{equation}
where $\Gamma(n)$ is the space of all polynomials of degree $\leq n$
in $N$ variables, and 
$\vol$ is the standard Euclidean volume on $\Gamma(n)\iso\C^{M(n)}$ with 
respect to its monomial basis (see \cite{Rumely:Lau}, Theorems 2.3 and 2.6, 
and the remarks on p.557).  

Suppose $k$ is a number field and $v$ is a place of $k$.   
Given a positive multiple $mH_0$ of the hyperplane at infinity of 
$(\C_v)^N$ and a bounded set 
$E_v \subset (\C_v)^N$ which is stable under 
the group of continuous automorphisms $\Gal_c(\C_v/k_v)$, we define the
{\em local sectional capacity}  $S_{\gamma}(E_v,mH_0)_v$ by 
a limit analogous to (\ref{S}):
\begin{equation} \label{SS}
\log S_\gamma(E_v,mH_0)_v = 
  -\lim_{n\to\infty} \frac{(N+1)!}{n^{N+1}} 
\log \vol_v \{f\in\Gamma_v(mn):  \|f\|_{E,v} \leq 1\},
\end{equation}
where $\Gamma_v(mn) = k_v\otimes_k \Gamma(mn)$ and $\vol_v$ is 
induced by the Haar measure on $k_v$ and the 
monomial basis of $\Gamma_v(mn)$.  
The limit was shown to exist in \cite{Rumely:Lau} in this setting 
and in \cite{Rumely:Lau:Varley} on algebraic varieties.   
It should be observed that if $E_v = D_v(0,1)$ is the unit polydisc in $(\C_v)^N$, 
then $S_\gamma(E_v,mH_0)_v=1$ \cite[p.555]{Rumely:Lau}.   

When $m = 1$, by \cite{Rumely:Lau}, Theorems 2.3 and 2.6, we have 
$S_{\gamma}(E_v,H_0) = d_{\infty}(E_v,H_0)$.  Hence by (\ref{DRel}),  
\begin{equation} \label{vpRel} 
S_{\gamma}(E_v,H_0)_v = d_{\infty}(E_v)_p^{N [k_v:\Q_p]} .
\end{equation} 
When $v$ is an archimedean place of a number field $k$, 
by \cite[p.557]{Rumely:Lau} 
\begin{equation} \label{Sgarch}
S_\gamma(E_v,H_0)_v = d_\infty(E_v)^{N[k_v:\R]} \ . 
\end{equation} 

\bigskip
\noindent{\bf Global sectional capacity.}
The sectional capacity has a global formulation for adelic sets of 
the form 
  $$\E = \prod_v E_v \subset \prod_v (\C_v)^N,$$
where the product is taken over all places $v$ of the global field $k$.  The 
global sectional capacity of $\E$ relative to the divisor $mH_0$ 
is the quantity $S_\gamma(\E, mH_0)$ defined by 
  $$\log S_\gamma(\E, mH_0) = \sum_v \log S_\gamma(E_v,mH_0)_v,$$
when the sum exists.  It is known to exist very generally 
\cite{Rumely:Lau:Varley}, but we will need
only the case when $E_v\subset (\C_v)^N$ is bounded and stable under 
$\Gal_c(\C_v/k_v)$ for all $v$ and $E_v = D_v(0,1)$
is the unit polydisc for all but finitely many $v$  \cite{Rumely:Lau}.  
 
 \begin{prop} \label{global}
Let $k$ be a number field and $F:\C^N\to\C^N$ a regular polynomial
map of degree $d$ with coefficients in $k$.  Let 
$\E = \prod_v E_v \subset \prod_v (\C_v)^N$ be a bounded adelic 
set such that $E_v$ is stable under $\Gal_c(\C_v/k_v)$ for all $v$ 
and $E_v = D_v(0,1)$ is the unit polydisc for all but finitely many $v$.  
Then,
 $$S_\gamma(F^{-1}\E, H_0) = S_\gamma(\E, H_0)^{1/d}.$$
\end{prop}

\proof  
From \cite{Rumely:Lau:Varley}, Theorem 10.1,  it follows that 
$S_\gamma(F^{-1}\E, F^*H_0) = S_\gamma(\E, H_0)^{d^N}$.
On the other hand, $S_\gamma(F^{-1}\E, F^*H_0) = 
S_\gamma(F^{-1}\E,d H_0) = S_\gamma(F^{-1}\E, H_0)^{d^{N+1}}$
by the homogeneity property of the sectional capacity, implicit in (\ref{SS}).  
\qed


\bigskip
\section{The pullback formula (over $\C$)}

In this section, we prove Theorem \ref{pullback}, 
the formula for the transfinite diameter of the preimage of a bounded set in $\C^N$
by a regular polynomial endomorphism.  We apply the 
pullback formula for the global sectional capacity (Proposition \ref{global})
and Theorem \ref{Julia}.

\begin{prop} \label{Qiconstant}
Let $F:\C^N\to\C^N$ be a homogeneous regular polynomial map of degree $d \ge 2$
with coefficients in $k=\Q(i)$.  Then
  $$d_\infty(F^{-1}E) = |\Res(F)|^{-1/Nd^N} d_\infty(E)^{1/d}$$
for all bounded sets $E\subset\C^N$.  
\end{prop}

\proof
Let $E\subset\C^N$ be a bounded set.  Define 
 $${\E} := E\times\prod_{v\not=\infty} D_v(0,1) \subset 
   \prod_v \A^N(\C_v)\subset \prod_v \P^N(\C_v).$$
The global sectional capacity of $\E$ (relative to $k = \Q(i)$
and the hyperplane at $\infty$, $H_0$) is then
 $$S_\gamma(\E,H_0) = S_\gamma(E,H_0)_\infty \times
 	\prod_{v\not=\infty} S_\gamma(D_v(0,1),H_0)_v
	= d_\infty(E)^{2N},$$
because $S_\gamma(D_v(0,1),H_0)_v = 1$ for all nonarchimedean $v$ of $k$
and $[k_{\infty}:\R]=2$.  

By the global pullback formula
of Proposition \ref{global}, 
  $$S_\gamma(F^{-1}\E, H_0) = S_\gamma(\E, H_0)^{1/d}.$$
Applying the local decompositions of both sides, we find that
  $$S_\gamma(F^{-1}E,H_0)_\infty \cdot \prod_{v\not=\infty}
  S_\gamma(F^{-1}D_v(0,1), H_0) = S_\gamma(E,H_0)_\infty^{1/d}.$$
Therefore, 
 $$d_\infty(F^{-1}E) = C_F  \cdot d_\infty(E)^{1/d}$$
with $C_F = \prod_{v\not=\infty} S_\gamma(F^{-1}D_v(0,1),H_0)_v^{-1/2N}$,
independent of the set $E$.  

On the other hand, for the filled Julia set $K_F$ of $F$, the
invariance $F^{-1}K_F = K_F$ implies that 
  $$C_F = d_\infty(K_F)^{(d-1)/d}.$$
By Theorem \ref{Julia}, we have $d_\infty(K_F) = |\Res(F)|^{-1/Nd^{N-1}(d-1)}$,
so $C_F = |\Res(F)|^{-1/Nd^N}$.   
\qed

\bigskip
\noindent{\bf Proof of Theorem \ref{pullback}.}
Let $F:\C^N\to\C^N$ be a regular polynomial endomorphism of degree $d$.  
When $d = 1$, Theorem \ref{pullback}
is equivalent to She\u{i}nov's formula \cite{Sheinov}, 
so we can assume that $d \ge 2$. 

Fix a bounded set $E \subset \C^N$ and 
let $E_\eps$ be the closed $\eps$-neighborhood of $E$.  Then $E_\eps$ is compact.
By the outer regularity of the transfinite diameter, $d_\infty(E_\eps)$ descends
to $d_\infty(E)$ as $\eps\to 0$ (\cite{Rumely:Lau}, Theorem 2.9).
As observed by Bloom and Calvi in \cite{Bloom:Calvi}, Theorem 5,
the transfinite diameter of a preimage depends only on the 
leading homogeneous part of $F$; that is, 
\begin{equation} \label{BCformula} 
          d_\infty(F^{-1}E) = d_\infty(F_h^{-1}E), 
\end{equation} 
and we can therefore assume that $F$ is homogeneous. (Although Bloom and Calvi
state their result for compact sets, if $\overline{E}$ is the closure of $E$,
trivially $d_{\infty}(\overline{E}) = d_{\infty}(E)$, 
so (\ref{BCformula}) holds for bounded sets.)  

Choose a sequence of homogeneous polynomial maps $F_n$ with 
coefficients in $\Q(i)$ such that the coefficients of $F_n$ converge
to the coefficients of $F$ as $n\to \infty$.  Then in fact, $F_n\to F$ 
uniformly on compact sets.  

For a fixed $\eps>0$, we have
  $$F_n^{-1}E_\eps \supset F^{-1}E$$
for all sufficiently large $n$.  Consequently, 
  $$d_\infty(F_n^{-1}E_\eps) \geq d_\infty(F^{-1}E).$$
From Proposition \ref{Qiconstant}, we know that the left hand side
equals $|\Res(F_n)|^{-1/Nd^N}d_\infty(E_\eps)^{1/d}$.  Letting
$n\to\infty$, the continuity of the resultant implies that 
  $$|\Res(F)|^{-1/Nd^N} d_\infty(E_\eps)^{1/d} \geq d_\infty(F^{-1}E).$$
As $\eps$ was arbitrary, we find that 
  $$|\Res(F)|^{-1/Nd^N} d_\infty(E)^{1/d} \geq d_\infty(F^{-1}E).$$

For the reverse inequality, again fix $\eps>0$.  
For all sufficiently large $n$, we have 
  $$F_n^{-1}E \subset (F^{-1}E)_\eps,$$ 
and consequently, 
 $$|\Res(F_n)|^{-1/Nd^N}d_\infty(E)^{1/d} 
   = d_\infty(F_n^{-1}E) \leq d_\infty( (F^{-1}E)_\eps).$$
Taking the limit as $n\to\infty$ and then letting $\eps\to 0$
completes the proof.  
\qed


\bigskip
\section{The pullback formula at nonarchimedean places}
\label{NonarchCase} 

In this final section, we prove Theorem \ref{pullbackv}.  We first need
two lemmas. 

\begin{lemma} \label{Res1}
Let $F : \C_p^N \to \C_p^N$ be a regular polynomial 
map of degree $d$. Suppose $F$ has coefficients in $\hO_p$, and that
$|\Res(F_h)|_p = 1$.  Then, writing $D_p(0,1)$ for the unit polydisc in $\C_p^N$, 
$$F^{-1}D_p(0,1) = D_p(0,1).$$
\end{lemma} 

\proof  Since $F$ has coefficients in $\hO_p$, 
$D_p(0,1) \subset F^{-1}D_p(0,1)$.  Since $|\Res(F_h)|_p = 1$, the induced map 
$\overline{F} = F \pmod{m_p}$ on $\tF_p^N$ is regular, which means that it extends to
a rational endomorphism of $\P^N(\tF_p)$ taking the hyperplane at infinity
to itself.  This, in turn, implies that $F$ takes $\P^N(\C_p) \backslash D_p(0,1)$
to itself, so $F^{-1}D_p(0,1)  = D_p(0,1)$.    
\qed 

\vskip  .1in 
If $k_v \subset \C_p$ is a local field and 
$F = \sum_{\alpha} c_{\alpha} z^{\alpha} \in k_v[z_1,\ldots,z_N]$, 
put $|F|_v = \max_{\alpha} |c_{\alpha}|_v$.   
Recall that if $L/k$ is a
finite extension of number fields, 
a place $v$ of $k$ is said to split completely
in $L$ if for every place $w$ of $L$ lying over $v$, we have $L_w \cong k_v$.  

\begin{lemma} \label{MBapprox}
Let $k_v \subset \C_p$ be a nonarchimedean local field, 
the completion of a number field $k$ at a place $v$.  
Suppose $F : \C_p^N \rightarrow \C_p^N$ is a regular polynomial map 
with coefficients in $k_v$. Then for any $\varepsilon > 0$, 
there is a finite extension $L/k$ in which $v$ splits completely, and
a regular polynomial map $\hF$ with coefficients in $L$,
such that 
\begin{enumerate} 
\item For each place $w$ of $L$ lying over $v$, we have
 $|\hF-F|_w < \varepsilon$ and $|\Res(\hF_h)|_w = |\Res(F_h)|_v$.   
\item For each nonarchimedean place $w$ of $L$ not lying over $v$, we have 
 $|\hF|_w \le 1$ and $|\Res(\hF_h)|_w = 1$.  
\end{enumerate} 
\end{lemma} 

\proof
Suppose $F$ has degree $d$.  
Identify each $N$-vector $G$ of polynomials of degree $d$ in $N$ variables
with its vector of coefficients, 
viewed as an element of an affine space $\A^M$; let $t$ be an additional variable,
and consider the subvariety
$V$ of $\A^{M+1}$ defined by the equation 
$$\Res(G_h) \cdot t = 1.$$
Then points of $V$ correspond to regular polynomial maps of degree $d$.
By \cite[p.252]{GKZ}, $V$ is an absolutely irreducible affine variety over $\Q$,  
defined by an equation with integer coefficients having no common divisor.  
Thus $V$ is the generic fibre of an irreducible 
scheme $\cV/\Spec(\Z)$, which surjects onto $\Spec(\Z)$.  
Put $V_k = V \times_\Q \Spec(k)$.
For each nonarchimedean place $u$ of $k$, 
write $\hO_u$ for the ring of integers of $\C_u$.  Then $V_k(\hO_u)$ is nonempty,
since it contains the point corresponding to $G_0(z) = (z_1^d,\ldots,z_N^d)$
with $\Res(G_0) = 1$.  
Note that each $G \in V_k(\hO_u)$ satisfies $|\Res(G_h)|_u = 1$. 
At the place $v$, the polynomial $F$ corresponds to a point of $V_k(k_v)$. 
Let $\Omega_v$ be the $\varepsilon$-neighborhood of $F$ in $V_k(k_v)$.  After
shrinking $\varepsilon$, if necessary, we can assume that all the polynomials $G$
corresponding to points in $\Omega_v$ satisfy $|\Res(G_h)|_v = |\Res(F_h)|_v$.  
   
By Moret-Bailly's ``Existence theorem for incomplete Skolem problems'' 
(\cite{MB}, Theorem 1.3), there is a finite extension $L/k$ in which $v$ 
splits completely, and a point $x \in V_k(L)$ 
whose $\Gal(L/k)$-conjugates (viewed as embedded in $V_k(\C_u)$ for each place
$u$ of $k$) belong to $V_k(\hO_u)$ for each nonarchimedean $u \ne v$ 
and belong to $\Omega_v$ if $u = v$. 
Let $\hF$ be the map corresponding to the point $x$.
\qed    


\bigskip\noindent {\bf Proof of Theorem \ref{pullbackv}.}
The proof has two steps.  First, we prove the formula for regular
maps $F: \C_p^N \rightarrow \C_p^N$ defined over a local field $k_v$ 
and for bounded sets
of the form $U = \{ z \in \C_p^N : |f_1(z)|_p \le 1, \ldots, |f_M(z)|_p \le 1 \}$ 
with $f_1, \ldots, f_M \in k_v[z]$.  Then, we use approximation properties of 
the transfinite diameter to deal with the general case.     

First suppose $F$ is defined over a local field $k_v \subset \C_p$, and that 
$U$ is bounded and of the form (\ref{Nbhd}) with the $f_i$ defined over $k_v$.  
Clearly $U$ is invariant under $\Gal_c(\C_p/k_v)$.  Put $W = F^{-1}U$; 
then $W$ is bounded (say $W \subset D_p(0,R)$) and is defined by the equations 
$|f_i \circ F(z)|_p \le 1$, $i = 1,\ldots,M$.  By continuity, 
there is an $\varepsilon > 0$ such that if $G$ is any regular polynomial map 
of degree $d$ defined over $k_v$ with $|G-F|_v < \varepsilon$, 
then $G^{-1}U \subset D_p(0,R)$, $|\Res(G_h)|_p = |\Res(F_h)|_p$, 
and $|f_i \circ F(z) - f_i \circ G(z)|_v < 1/2$ for all $z \in D_p(0,R)$, 
for each $i$.  By the ultrametric inequality, it follows that 
$$ G^{-1}U = F^{-1}U.$$

Let $k$ be a number field with a place $v$ inducing
$k_v$. By Proposition \ref{MBapprox}, there is a finite extension $L/k$ in 
which $v$ splits completely, and a map $\hF$ defined over $L$, such that 
for each place $w$ of $L$ 
\begin{enumerate}
\item if $w$ lies over $v$, then using the 
isomorphism $L_w \cong k_v$ to identify $\hF$ with a map defined over $k_v$, 
we have $\hF^{-1}(U) = F^{-1}U$ and $|\Res(\hF_h)|_w = |\Res(F_h)|_v$;
\item for each nonarchimedean  $w$ not over $v$, then $|\hF|_w \le 1$ 
and $|\Res(\hF_h)|_w = 1$.  
\end{enumerate}  

Define an adelic set $\E_L = \prod_w E_w$ by taking
$E_w = U$ for each $w$ lying over $v$,  $E_w = D_w(0,1) \subset \C_w^N$ 
if $w$ is nonarchimedean and does not lie over $v$, 
and let $E_w \subset \C^N$ be a set with 
$S_{\gamma}(E_w,H_0)_w = 1$ (for example, the unit polydisc;  
see \cite{Rumely:Lau}, Example 4.3, p.558) if $w$ is archimedean.  
Then $S_{\gamma}(E_w,H_0)_w = 1$ for each $w$ not over $v$.
By Lemma \ref{Res1}, if $w$ is nonarchimedean 
and does not lie over $v$, 
$$
S_{\gamma}(\hF^{-1}E_w,H_0)_w = 1.
$$
If $w$ is archimedean, by Theorem \ref{pullback} and formula (\ref{Sgarch}), 
$$
S_{\gamma}(\hF^{-1}E_w,H_0)_w = |\Res(\hF_h)|_w^{-1/d^N}.
$$

There are $[L:K]$ places of $L$ over $v$, 
and $S_{\gamma}(\hF^{-1}U,H_0)_w = S_{\gamma}(F^{-1}U,H_0)_v$ 
for each $w$ over $v$.       
Applying the pullback formula for the global sectional capacity 
(Proposition \ref{global}) and using the local decomposition of each side, we get
$$
 (\prod_{w|\infty} |\Res(\hF_h)|_w^{-1/d^N})
\cdot S_{\gamma}(F^{-1}U,H_0)_v^{[L:K]} = 
(S_{\gamma}(U,H_0)_v^{1/d})^{[L:K]}.
$$
By the product formula and the fact that $|\Res(\hF_h)|_w = 1$ 
for all nonarchimedean $w$ not over $v$, while 
$|\Res(\hF_h)|_w = |\Res(F_h)|_v$ for each $w$ over $v$, 
$$
(\prod_{w|\infty} |\Res(\hF_h)|_w)^{-1} = (|\Res(F_h)|_v)^{[L:K]}.
$$
Combining the last two formulas gives 
\begin{equation} \label{Sfp} 
S_{\gamma}(F^{-1}U,H_0)_v = |\Res(F_h)|_v^{-1/d^N} 
\cdot S_{\gamma}(U,H_0)_v^{1/d}. 
\end{equation}
However, $|\Res(F_h)|_v = |\Res(F_h)|_p^{[k_v:\Q_p]}$.
Thus by (\ref{vpRel}), formula (\ref{Sfp}) is equivalent to 
\begin{equation} \label{Dfp} 
d_{\infty}(F^{-1}U)_p = |\Res(F_h)|_p^{-1/Nd^N} \cdot d_{\infty}(U)_p^{1/d}  
\end{equation} 
as was to be shown.

\vskip .1 in
Now let $F : \C_p^N \rightarrow \C_p^N$ be an arbitrary regular polynomial 
map of degree $d$.  If $U$ is a set of the form (\ref{Nbhd}), 
then since $\tQ_p$ is dense in $\C_p$, we can assume that the functions $f_i$
determining $U$ are defined over $\tQ_p$.  For the same reason, 
there is a map $\hF$ defined over $\tQ_p$ with $\hF^{-1}U = F^{-1}U$ 
and $|\Res(\hF_h)|_p = |\Res(F_h)|_p$.  Since these functions have only 
finitely many coefficients, they are actually defined over a local field $k_v$.  
By (\ref{Dfp}), 
$$
d_{\infty}(F^{-1}U)_p = d_{\infty}(\hF^{-1}U)_p 
= |\Res(F_h)|_p^{-1/Nd^N} \cdot d_{\infty}(U)_p \ . 
$$
Thus, (\ref{Dfp}) holds without restriction on $F$ and $U$.  

Let an arbitrary bounded set $E \subset \C_p^N$ be given.  
For each $\varepsilon > 0$, there is a neighborhood $U$ of $E$ 
of the form (\ref{Nbhd}) with $d_{\infty}(U)_p \le d_{\infty}(E)_p + \varepsilon$.  
Applying (\ref{Dfp}) and letting $\varepsilon \rightarrow 0$ gives 
\begin{equation} \label{UpIneqp}
d_{\infty}(F^{-1}E)_p \le |\Res(F_h)|_p^{-1/Nd^N} \cdot d_{\infty}(E)_p \ .
\end{equation} 

For the reverse inequality, apply the approximation property to $F^{-1}E$.  
For each $\varepsilon > 0$, there is a neighborhood $U$ of $F^{-1}E$ of
the form (\ref{Nbhd}) such that 
$d_{\infty}(U)_p \le d_{\infty}(F^{-1}E)_p + \varepsilon$.  
We will construct a polynomially-defined 
neighborhood $V$ of $E$ with $F^{-1}V \subset U$ 
by first pulling $U$ back to a galois cover, then intersecting $U$ with its
conjugates to get a galois-invariant set $V_X$, 
and finally descending $V_X$ to obtain $V$.   

Writing $Y = Z = \P^N$, 
extend $F$ to a map of normal varieties $F : Y \rightarrow Z$ over $\C_p$, 
and let $F^* : \C_p(Z) \rightarrow \C_p(Y)$  be the induced map on the function
fields.  Let $\cK$ be the galois closure of $\C_p(Y)/\C_p(Z)$, and let $X$ be the 
normalization of $Z$  in $\cK$.  Then $\cK = \C_p(X)$.  
As $F$ is finite, the canonical map $X \rightarrow Z$ factors through a finite map 
$G : X \rightarrow Y$ such that the maps on points 
$$ 
 X(\C_p) \stackrel{G}{\longrightarrow} Y(\C_p) 
           \stackrel{F}{\longrightarrow} Z(\C_p) 
$$
correspond contravariantly to the field inclusions  
$\C_p(Z) \hookrightarrow \C_p(Y) \hookrightarrow \C_p(X)$.  
Put $E_X = (F \circ G)^{-1}E$.  
The galois group $\cG = \Gal(\C_p(X)/\C_p(Z))$ acts on $X(\C_p)$, 
and stabilizes $E_X$.
 
Let $f_1, \ldots, f_M \in \C_p(Y)$ be the functions determining $U$;  
by abuse of notation (identifying $f_i$ with $f_i \circ G$) 
we can view them as elements of $\C_p(X)$.  Let
\begin{equation*}
V_X = \{x \in X(\C_p) : |\sigma(f_i)(x)|_p \le 1  
\text{ for all $i = 1, \ldots, M$ and all $\sigma \in \cG$} \} .
\end{equation*} 
Then $V_X$ is a neighborhood of $E_X$ contained in $G^{-1}U$, stable
under $\cG$.  For each $f_i$ (viewed as an element of $\C_p(Y)$), let
$$
P_i(t) = t^{D_i} + a_{i,1}(z) t^{D_i-1} + \cdots + a_{i,D_i}(z) 
$$
be its minimal polynomial over $\C_p(Z)$.  
Then the $a_{i,j}(z)$ are elementary symmetric functions of the $\sigma(f_i)$
which are polynomials in $z_1, \ldots, z_N$. 
Since $E_X \subset V_X$, the ultrametric inequality shows 
that $|a_{i,j} \circ (F \circ G)(x)|_p \le 1$ for each $x \in E_X$. 
But $F \circ G$ maps $E_X$ onto $E$, 
so $|a_{i,j}(z)|_p \le 1$ for each $z \in E$.  
Hence
\begin{equation*} 
V := \{ z \in Z(\C_p) : |a_{i,j}(z)|_p \le 1 \text{ for all $i, j$} \} 
\end{equation*}
is a neighborhood of $E$ in $\C_p^N$.  

On the other hand, the theory of Newton polygons (see \cite{Artin}, \S 2.5) 
shows that for each fixed $z \in V$, the roots $\alpha_{i,\ell}$ 
of $P_i(t) = P_{i,z}(t)$ satisfy $|\alpha_{i,\ell}|_p \le 1$.  
For each $x \in X(\C_p)$ with $F \circ G(x) = z$,
the $\alpha_{i,\ell}$ are precisely the values 
$\sigma(f_i)(x)$ for $\sigma \in \cG$.   
It follows that 
$$(F \circ G)^{-1}V = V_X \subset  G^{-1}U\ ,$$ 
so $F^{-1}V \subset U$.    
This gives 
\begin{eqnarray*}
|\Res(F_h)|_p^{-1/Nd^N} d_{\infty}(E)_p 
&\le& |\Res(F_h)|_p^{-1/Nd^N} d_{\infty}(V)_p\\
&=& d_{\infty}(F^{-1}V)_p \ \le \ d_{\infty}(F^{-1}E)_p + \varepsilon \ ,  
\end{eqnarray*}
and since $\epsilon > 0$ is arbitrary, we are done.  
\qed

\bigskip

\begin{thebibliography}{GKZ}

\bibitem[Ah]{Ahlfors:conformal}
L.~Ahlfors.
\newblock {\em Conformal {I}nvariants: {T}opics in {G}eometric {F}unction
  {T}heory}.
\newblock McGraw-Hill Book Co., New York, 1973.

\bibitem[Ar]{Artin}
E.~Artin.
\newblock {\em Algebraic functions and algebraic numbers}. 
\newblock Gordon and Breach, New York, 1967.

\bibitem[BR]{Baker:Rumely:equidistribution}
M.~Baker and R.~Rumely.
\newblock {Equidistribution of small points, rational dynamics, and potential
  theory}.
\newblock {To appear in \em Annales de l'Institut Fourier}.

\bibitem[BB]{Bassanelli:Berteloot}
G.~Bassanelli and F.~Berteloot.
\newblock {Bifurcation currents in holomorphic dynamics on ${\bf P}^k$, 2005}.
\newblock {Preprint available at \verb+http://arxiv.org/abs/math.DS/0507555+}. 

\bibitem[BJ]{Bedford:Jonsson}
E.~Bedford and M.~Jonsson.
\newblock {Dynamics of regular polynomial endomorphisms of ${\bf {C}}\sp k$}.
\newblock {\em Amer. J. Math.} {\bf 122} (2000), 153--212.

\bibitem[BT]{Bedford:Taylor}
E.~Bedford and B.A.~Taylor.
\newblock {A new capacity for plurisubharmonic functions}.
\newblock {\em Acta Math.} {\bf 149} (1982), 1--40.

\bibitem[BC]{Bloom:Calvi}
T.~Bloom and J.-P. Calvi.
\newblock {On the multivariate transfinite diameter}.
\newblock {\em Ann. Polon. Math.} {\bf 72} (1999), 285--305.

\bibitem[BD]{Briend:Duval}
J.-Y. Briend and J.~Duval.
\newblock {Exposants de {L}iapounoff et distribution des points p\'eriodiques
  d'un endomorphisme de $\bf {C}{\rm {P}}\sp k$}.
\newblock {\em Acta Math.} {\bf 182} (1999), 143--157.

\bibitem[De]{Demailly:lelong}
J.-P.~Demailly.
\newblock {Monge-{A}mp\`ere operators, {L}elong numbers and intersection
              theory}.
\newblock In {\em Complex analysis and geometry}, pages 115--193.
  Univ. Ser. Math.  Plenum, New York, 1993.

\bibitem[DeM]{D:lyap}
L.~DeMarco.
\newblock {Dynamics of rational maps: {L}yapunov exponents, bifurcations, and
  capacity}.
\newblock {\em Math. Ann.} {\bf 326} (2003), 43--73.

\bibitem[FS]{Fornaess:Sibony}
J.~E. Forn{\ae}ss and N.~Sibony.
\newblock {Complex dynamics in higher dimensions}.
\newblock In {\em Complex Potential Theory (Montreal, PQ, 1993)}, pages
  131--186. Kluwer Acad. Publ., Dordrecht, 1994.

\bibitem[GKZ]{GKZ}
I.M. Gel{\cprime}fand, M.M. Kapranov, and A.V. Zelevinsky.
\newblock {\em Discriminants, resultants, and multidimensional determinants}.
\newblock Mathematics: Theory \& Applications. Birkh\"auser Boston Inc.,
  Boston, MA, 1994.

\bibitem[Go]{Goluzin}
G.~M. Goluzin.
\newblock {\em Geometric theory of functions of a complex variable}.
\newblock Translations of Mathematical Monographs, Vol. 26. American
  Mathematical Society, Providence, R.I., 1969.

\bibitem[HP]{Hubbard:Papadopol}
J.~Hubbard and P.~Papadopol.
\newblock {Superattractive fixed points in ${\bf C}^n$}.
\newblock {\em Indiana Univ. Math. J.} {\bf 43} (1994), 321--365.

\bibitem[Je]{Jedrz}
M.~J{\polhk{e}}drzejowski.
\newblock {The homogeneous transfinite diameter of a compact subset of
              {${\bf C}\sp N$}}.
\newblock In {\em Proceedings of the Tenth Conference on Analytic Functions
              (Szczyrk, 1990)}.
\newblock   {Ann. Polon. Math.} {\bf 55} (1991), 191--205.

\bibitem[Kl]{Klimek}
M.~Klimek.
\newblock {\em Pluripotential Theory}.
\newblock The Clarendon Press Oxford University Press, New York, 1991.
\newblock Oxford Science Publications.

\bibitem[MB]{MB}
L.~Moret-Bailly.
\newblock {Groupes de Picard et probl\`emes de Skolem II}.
\newblock {\em Ann scient. \'Ec. Norm. Sup.} 22 (1989), 181--194.

\bibitem[Ra]{Ransford}
T.~Ransford.
\newblock {\em Potential theory in the complex plane}.
\newblock{Cambridge University Press, Cambridge, 1995} 

\bibitem[Ru]{Rumely:Vformula}
R.~Rumely.
\newblock {A Robin formula for the Fekete-Leja transfinite diameter, 2004}.
\newblock {Preprint available at \verb+http://arxiv.org/abs/math.CO/0407427+}.  

\bibitem[RL]{Rumely:Lau}
R.~Rumely and C.~F. Lau.
\newblock {Arithmetic capacities on {$\bf P\sp N$}}.
\newblock {\em Math. Z.} {\bf 215} (1994), 533--560.

\bibitem[RLV]{Rumely:Lau:Varley}
R.~Rumely, C.~F. Lau, and R.~Varley.
\newblock {Existence of the sectional capacity}.
\newblock {\em Mem. Amer. Math. Soc.} {\bf 145} (2000).

\bibitem[Sh]{Sheinov}
V.~P. She{\u\i}nov.
\newblock {An invariant form of {P}\'olya's inequalities}.
\newblock {\em Siberian Math. J.} {\bf 14} (1973), 138--145.


\bibitem[VdW]{vdW}
B.~L. van der Waerden.
\newblock {\em Modern Algebra, 2nd edition, vol. 2 (tr. Benak)}.
\newblock Ungar, New York, 1950.

\bibitem[Za]{Zaharjuta}
V.~P. Zaharjuta.
\newblock {Transfinite diameter, \v {C}eby\v sev constants and capacity for a
  compactum in ${\bf C}^n$}.
\newblock {\em Mat. Sb. (N.S.)} {\bf 96(138)} (1975), 374--389, 503.

\end{thebibliography}
\def\cprime{$'$}
\def\polhk#1{\setbox0=\hbox{#1}{\ooalign{\hidewidth
    \lower1.5ex\hbox{`}\hidewidth\crcr\unhbox0}}}

\end{document}